\documentclass[11pt,twoside]{article}
\usepackage{latexsym,amssymb}
\usepackage{rotating}

\date{\today}

\author{Peter Ashwin\\
School of Mathematical Sciences,\\
Laver Building, \\
University of Exeter, \\
Exeter EX4 4QE, UK.
\and
James Montaldi\\
CNRS-INLN\\
1361 route des Lucioles\\
06560 Valbonne,
France}

\title{Group theoretic conditions for existence of robust relative
  homoclinic trajectories}

\newtheorem{theorem}{Theorem}[section]

\newtheorem{proposition}[theorem]{Proposition}
\newtheorem{corollary}[theorem]{Corollary}
\newtheorem{definition1}[theorem]{Definition}
\newtheorem{remark1}[theorem]{Remark}
\newtheorem{example1}[theorem]{Example}
\newenvironment{proof}%
{\addvspace\baselineskip\noindent\textsc{Proof}\quad}%
{\hspace*{\fill}$\Box$\par\addvspace\baselineskip}

\newenvironment{remark}{\begin{remark1}\rm}{\end{remark1}}
\newenvironment{example}{\begin{example1}\rm}{\end{example1}}

\newcommand{\C}{\mathbf{C}}
\newcommand{\D}{\mathbf{D}}
\newcommand{\R}{\mathbf{R}}
\newcommand{\Z}{\mathbf{Z}}

\newcommand{\vbar}{\overline{v}}

\newcommand{\gammabar}{\overline{\gamma}}
\newcommand{\Fix}{\mathop{\rm Fix}\nolimits}

\newcommand{\RHT}{\textsc{rht}}
\newcommand{\rht}{\textsc{rht}}

\newcommand{\Kginv}{K^{g^{-1}}}

\begin{document}
\maketitle

\begin{abstract}
  In this paper we consider robust relative homoclinic trajectories (RHTs)
  for $G$-equivariant vector fields. We give some
  conditions on the group and representation that imply existence of 
  equivariant vector fields with such 
  trajectories. Using these result we show very simply that abelian
  groups cannot exhibit relative homoclinic trajectories.
  Examining a set of group theoretic conditions that imply existence of
  RHTs, we construct some new examples of robust
  relative homoclinic trajectories. We also classify 
  RHTs of the dihedral and low order symmetric groups by means of their
  symmetries.
\end{abstract}

%%%%%%%%%%%%%
\section{Introduction}

Although homoclinic and heteroclinic cycles are not structurally
stable to perturbations in typical systems, they do play a very
important role in organising the dynamics of trajectories that are
nearby in phase and parameter space.  In the presence of symmetries
(or conserved quantities) it has been observed for many years that
extra structure can force homoclinic orbits to be robust (i.e.\ they
can persist under an open set of perturbations) and this gives their
study in equivariant systems a special importance.

More precisely, suppose we have a flow generated by 
$$
\dot{x}=v(x)
$$
where $x\in M$ is a smooth manifold where a compact Lie group $G$
acts smoothly is such that the flow generate by $v$ commutes with the
action of $G$. We say $X_0=G\cdot x_0$ is a {\em relative equilibrium}
if it is flow invariant, i.e.\ if $v(x_0)\in T (G\cdot x_0)$. A {\em
  relative homoclinic trajectory} (\RHT) is a trajectory $\gamma(t)$
with $\gamma(t)\not\in X_0$ but such that $\gamma(t)\rightarrow X_0$
as $t\rightarrow \pm\infty$.

In the case of a finite group $G$ a relative equilibrium consists of a
finite number of equilibrium points related to one another by the
symmetry group, and a homoclinic trajectory consists of a number of
trajectories that are homoclinic or more generally heteroclinic
between different points in the same relative equilibrium. More
generally if $G$ is continuous it will be a (usually infinite) number
of connections between equilibria or invariant tori in $X_0$. We say a
homoclinic trajectory is {\em robust} if it persists under an open set
of sufficiently smooth perturbations.

Most studies of robust \RHT s in equivariant systems have been either
general results that assume the existence of certain structures (for
example the conditions in \cite{KM} that imply attractivity of cycles)
or they are specific examples where a vector field or normal form near
bifurcation is studied in detail.  In this paper we look at two
general questions that are intimately related:
\begin{itemize}
\item How can we characterise the geometry of a relative homoclinic 
  orbit?
\item How can we characterise group actions that {\em allow} structurally
  stable relative homoclinic trajectories?  
\end{itemize} 
After reviewing two motivating examples of Guckenheimer and Holmes
\cite{Guc&Hol88} and of Kevrekidis {\em et al} \cite{Kev&al90} in
the remainder of this introduction, in Section~\ref{secgeom} we
introduce notation and describe the geometry of an \RHT\ and a coset
in the group that we term the {\em twist} of a \RHT. This is used in
Proposition~\ref{proptwist} to prove that a necessary condition for
existence of a robust \RHT\ is that the isotropies of trajectory,
endpoints and twist are related in a certain way. As a
trivial consequence, we deduce that one cannot have robust \RHT s
for abelian group actions.

In Section~\ref{sec:classification} we investigate the
group theoretic conditions in Proposition~\ref{proptwist} necessary
for existence of robust \RHT s. To this end, we introduce the notion of a
{\em homoclinic triple}; this is a triple $(K,g,H)$ where $K$ and $H$
are subgroups and $g$ is a group element with certain properties
that characterise the symmetries of a robust \RHT. After demonstrating
some useful invariance properties of homoclinic triples we classify in 
Theorem~\ref{thm:D4n} the set of all possible homoclinic triples
in the dihedral groups $\D_n$. This is followed by a discussion of
homoclinic triples in the symmetric groups $S_n$ and also in some wreath
product examples.

Section~\ref{secconst} adds extra hypotheses such that sufficient (but
not necessary) conditions on a group action are found for robust \RHT
s; one of these conditions is group theoretic while the other concerns
the geometry of the group representation. This is used to construct a
new example of an action of $SO(3)$ on a 16-dimensional vector space
with robust \RHT s and also a large class of finite groups that admit
actions with robust \RHT s. These include many with wreath product
structure \cite{Dia&al99}. The section ends with an example that
explores the gap between the necessary and sufficient conditions for
robust \RHT s. Finally, Section~\ref{secdiscuss} discusses extensions
and consequences of this work.

For the remainder of this introductory section we recall two standard
examples of systems and symmetry groups that permit robust relative
homoclinic trajectories, partly for motivation, but also for later
discussion.

\subsection{The Guckenheimer-Holmes robust {\sc RHT}}
\label{secGHmodel}

Consider the flow on $\R^3$ generated by
\begin{equation}\label{eqGHmodel}
  \begin{array}{rcl}
    \dot{x}_1 &=& x_1( \lambda-x_1^2+bx_2^2+cx_3^2)\\
    \dot{x}_2 &=& x_2( \lambda-x_2^2+bx_3^2+cx_1^2)\\
    \dot{x}_3 &=& x_3( \lambda-x_3^2+bx_1^2+cx_2^2)
  \end{array}
\end{equation}
where $b\neq c$ are real constants. This is equivariant under the
group $G$ generated by $\kappa:(x_1,x_2,x_3)\mapsto(-x_1,x_2,x_3)$ and
$\rho(x_1,x_2,x_3) \mapsto (x_2,x_3,x_1)$ (we could denote the group
$G=\Z_2(\kappa)\wr \Z_3(\rho)$ using the wreath product notation, see
\cite{Dio&al96}, or $\mathbb{T}_h$ in the Schoenflies notation).  For
certain open sets of $(b,c)$ with $bc<0$, Guckenheimer and Holmes
\cite{Guc&Hol88} observed and proved the existence of relative
equilibria $G\cdot(1,0,0)$ connected by a relative homoclinic
trajectory that is robust to any perturbation preserving the symmetry
$G$, and in particular to the addition of higher order polynomial
terms for the above vector field near $\lambda=0$. The equilibria in
question have symmetry $\Z_2(\kappa)\times\Z_2(\rho\kappa\rho^{-1})$
whereas the connections have symmetry $\Z_2(\kappa)$ (up to
conjugation).

\subsection{The Kuramoto-Sivashinsky robust {\sc rht}}\label{secKSmodel}

As noticed by Kevrikidis, Nicolaenko and Scovel \cite{Kev&al90} the
Kuramoto-Sivashinsky equation
$$
v_t+4v_{xxxx}+\alpha(v_{xx}-v v_x)=0
$$
posed on the line $x\in[0,2\pi]$ with periodic boundary conditions
can display attracting \RHT s for a certain range of parameter values.
This system commutes with the group $O(2)$ of symmetries generated by
translation $x\mapsto x+\theta$ modulo $2\pi$, and reflection
$x\mapsto 2\pi-x$. This equation exhibits robust \RHT s which can be
seen in numerical simulations. In an $N$-mode 
truncated Fourier representation this can be seen by writing
$$
v(x,t)=\sum_{n=1}^{N}\left(y_k(t) \cos kx+y_{N+k}(t) \sin kx\right)
$$
and then for $\alpha\in (16.13,22.557)$ one can observe (in the
truncated equations) \RHT s connecting relative equilibria with
$y_k=0$ for all $k=1,\cdots, N$ and $y_{N+k}=0$ for all odd
$k=1,3,5,\cdots, N$. This \RHT\ is simply the group orbit of a \RHT\ 
contained in the fixed point subspace of the subgroup generated by
$\kappa_1:x\rightarrow 2\pi-x$, $\kappa_2:x\rightarrow \pi-x$ and
$\rho:x\rightarrow x+\pi$ (all taken modulo $2\pi$ on the domain). In
this subspace the equilibria have symmetry $\Z_2(\kappa_1)\times
\Z_2(\kappa_2)$ and the connections have symmetry $\Z_2(\kappa_1)$
only. These \RHT s differ from the Guckenheimer-Holmes example in that
any connections between a pair of equilibria come in pairs.

%%%%%%%%%%%%%

\section{The geometry of relative homoclinic trajectories}
\label{secgeom}

We consider relative homoclinic trajectories with reference to the
flow on the orbit space $M/G$. Recall that the orbit space $M/G$ has a
natural stratification by orbit (isotropy) type, whose strata we
denote by $(M/G)_{(H)}$, where $(H)$ is the conjugacy class of the
subgroup $H$ of $G$.  

We now recall some standard notation we need and then state the
problem more precisely. The normalizer, $N_G(H)$, of $H$ in $G$ is the
group of $g\in G$ such that $gH=Hg$. $\Fix(H)$ is the fixed point
space of $H$ in $M$, i.e.\ the set of $x\in M$ such that $hx=x$ for
all $h\in H$.  If $H$ is an isotropy subgroup then $g$ maps $\Fix(H)$
to itself if and only if $g\in N_G(H)$ \cite{MDG93}.  For $g\in G$, we
write $H^g$ to mean the conjugate subgroup by $g$, i.e.\
$H^g=gHg^{-1}$. We also use the notation $H<G$ to mean that $H$ is a
subgroup of $G$.

Let $v$ be a $G$-equivariant vector field. The flow on $M$ descends to
a flow on the orbit space $M/G$, induced by a stratified vector field
$\vbar$.  A {\em relative homoclinic trajectory}, or \RHT, of the
dynamics on $M$ is a trajectory whose image in $M/G$ is simply a
homoclinic trajectory.  For an \RHT\ $\gamma(t)$ we denote its image
in $M/G$ by $\gammabar(t)$. The $\alpha$- and $\omega$-limit sets of
the image $\gammabar(t)$ are then the same (relative) equilibrium.

Let $\gamma(t)$ be an \RHT.  Then each point $\gamma(t)$ has the same
isotropy subgroup, which we call the isotropy subgroup of $\gamma$ and
denote by $K$.

The limit sets $\alpha(\gamma)$ and $\omega(\gamma)$ are closed,
connected flow-invariant subsets of the relative equilibrium
$\alpha(\gammabar) = \omega(\gammabar)$.  These flow-invariant subsets
also have an isotropy, and we denote the isotropy of $\alpha(\gamma)$
by $H$.  By continuity we have that $K<H$.

If the group $G$ is finite, then $\alpha(\gamma)$ is an equilibrium
point $x$, and $\omega(\gamma) = g\cdot x$, for some $g\in G$.  This
element $g$ we call the \emph{twist} of the \rht. More properly, it is
\emph{a twist}, for it is only well-defined modulo $H$; that is, the
twist is naturally an element of $G/H$, for if $y=g\cdot x$ then
$y=gh\cdot x$ for all $h\in H$.  The equilibrium point
$\omega(\gamma)$ is of isotropy $H^g=gHg^{-1}$, and by continuity
again $K<H^g$.

\begin{figure}
\begin{center}
  \begin{picture}(100,60)(-45,-10)
    %% Fix H:
    \put(0.5,10){\mbox{$\Fix(H)$}}
    \put(0,0){\line(0,1){50}}
    \put(40.5,10){\mbox{$\Fix(H^g)$}}

    %% Fix K
    \put(15,45){\mbox{$\Fix(K)$}}
    \put(0,0){\line(1,0){40}}
    \put(0,50){\line(1,0){40}}
    \put(40,0){\line(0,1){50}}

    %% Fix K'
    \put(-20,35){\mbox{\begin{rotate}{25}$\Fix(K^{g^{-1}})$\end{rotate}}}
%    \put(-27,-3){\mbox{$\Fix(K^{g^{-1}})$}}
    \put(-30,-15){\line(2,1){30}}
    \put(-30,-15){\line(0,1){50}}
    \put(-30,35){\line(2,1){30}}

    %% x etc
    \put(0,25){\circle*{1}}
    \put(-3,22){\mbox{$x$}}
    \put(40,25){\circle*{1}}
    \put(41,24){\mbox{$g\cdot x$}}
    \put(-30,10){\circle*{1}}
    \put(-42,9){\mbox{$g^{-1}\cdot x$}}

    %% gamma
    \qbezier(1,25.5)(20,34.5)(39,25.5)
    \put(20,33){\mbox{$\gamma$}}
    \put(18,29.4){\mbox{\large$\rightarrow$}}

    %% gamma'
    \qbezier[50](-29,11)(-15,30)(-1,25.5)
    \put(-18,17){\begin{rotate}{25}\mbox{$g^{-1}\gamma$}\end{rotate}}
    \put(-16,22){\begin{rotate}{25}\mbox{\large$\rightarrow$}\end{rotate}}

  \end{picture}
\end{center}

\caption{\label{figRHT}
  Schematic diagram showing (part of) a relative homoclinic trajectory
  to a relative equilibrium $x$ in the phase space $M$. Note that
  there may be many connections to and from $x$ in the \RHT; indeed
  there may also be many connections between $x$ and $g\cdot x$. In
  this case, the \RHT\ has isotropy $K$ and twist $g$.  }
\end{figure}
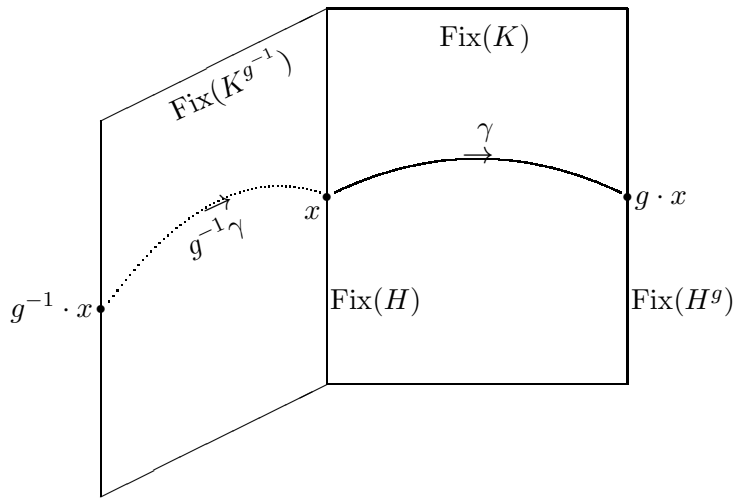

If on the other hand, $G$ is merely compact, $\alpha(\gamma)$ is of
the form $T\cdot x$, for any $x\in\alpha(\gamma)$, and some torus
$T<N_G(H)<G$ (not necessarily maximal).  Moreover, $g\cdot
x\in\omega(\gamma)$, for some $g\in G$, and $\omega(\gamma)=gT\cdot x
= (gTg^{-1})\cdot(g\cdot x)$.  In either case, we call $g$ the {\em
  twist} of the \RHT, and in this case it is well defined modulo
$(G/H)/T$.  Again the points of $\omega(\gamma)$ have isotropy $H^g$,
and again $K<H^g$. Figure~\ref{figRHT} schematically shows the setup
within $M$.

Thus, associated to any \RHT, there are subgroups $K$ and $H$, and a
twist $g\in G$ that is well-defined modulo $H$, or modulo $HT$ in the
compact case.  We write this triple as $(K,g,H)$.  If $\gamma$ is an
\rht\ with triple $(K,g,H)$ then for $f\in G$ the \rht\ $f\cdot
\gamma$ has triple $(K^f,g^f,H^f)$ as is readily checked.  Note that
$x$ will be a limit point for several distinct trajectories in
the \rht, and indeed
there may be more than one \rht\ between two given points $x$ and
$y=g\cdot x$, as in the Kuramoto-Sivashinsky example.

\begin{example}
  In the Guckenheimer-Holmes example in Section~\ref{secGHmodel}, the
  isotropy subgroup of the \RHT\ is $K=\Z_2(\kappa)$, while that of
  the equilibrium $x=(0,0,1)$ is $H\simeq\Z_2(\kappa) \times
  \Z_2(\rho\kappa\rho^{-1})$. The twist is then simply $g=\rho$.  For
  an open set of coefficients, the \RHT\ from $(0,0,1)$ tends to
  $(1,0,0)$, which is $g\cdot(0,0,1)$ for
  $g(x_1,x_2,x_3)=(x_3,x_1,x_2)$.  Notice that $g\not\in N_G(K)$.
\end{example}

\begin{example}
  For the Kuramoto-Sivashinsky example in Section~\ref{secKSmodel} we
  have $K=\Z_2(\kappa_1)$, $H=\Z_2(\kappa_1)\times\Z_2(\kappa_2)$ and
  $g=\rho$. Hence $K\neq K^g$ but $H=H^g$; in fact $H$ is in a
  conjugacy class of its own because its normalizer is the whole
  group.
\end{example}

We only consider homoclinic orbits in $M$ that are homoclinic to
normally hyperbolic relative equilibria, i.e.\ where the centre
manifold $W^c(x)$ is contained in the group orbit $G\cdot x$ and so
$\gamma\subset W^u(x)\cap W^s(g\cdot x)$. This is not a restriction
because for generic equivariant vector fields, relative equilibria are
normally hyperbolic (for example, see \cite[Lemma 3.1.10]{Fie96m}).

The following result gives a simple necessary condition for an \RHT\
to be robust.

\begin{proposition} \label{proptwist}
  Let $\gamma$ be a robust \RHT\  with isotropy type $K$ and twist $g$,
  and let $H$ be the isotropy of $\alpha(\gamma)$.  Then
  \begin{itemize}
  \item[(i)] $gH\cap N_G(K)=\emptyset$ and
  \item[(ii)] $K < H\cap H^g$.
  \end{itemize}
\end{proposition}

The first property implies both that $K\neq K^{g^{-1}}$ and that $H$
is a proper subgroup of $G$.  The second property is equivalent to
requiring that $H$ strictly contains both $K$ and $K^{g^{-1}}$.

\begin{remark}
  The quotient $(H\cap H^g)/K$ tells us about the number of
  heteroclinic connections between $x$ and $g\cdot x$ within a given
  \rht, for if $\alpha(\gamma) = \alpha(f\cdot\gamma)=x$ and
  $\omega(\gamma) = \omega(f\cdot\gamma) = g\cdot x$ then $f\in H\cap
  H^g$.  Moreover, $f\cdot\gamma = \gamma$ if and only if $f\in K$.
\end{remark}

\begin{proof} \emph{(of Proposition~\ref{proptwist})}
  For (i) it is enough to show that $g\not\in N_G(K)$, for $g$ is an
  arbitrary twist.  Suppose first that $g\in N_G(K)$ and so $g$ maps
  $\Fix(K)$ to itself. Then all equilibria and connections lie within
  $M'=\Fix(K)$ and so we consider the flow on $M'$, which is
  equivariant under the action of $N_G(K)/K$. Note that the action of
  $N_G(K)/K$ is free at all points $\gamma(t)$.
  
  Let $x\in\alpha(\gamma)$. Note that $g\cdot(W^s(x)\cap
  M')=W^s(g\cdot x) \cap M'$ and so in particular $\dim W^s(x)\cap
  M'=\dim W^s(g\cdot x)\cap M'$.  Since $\dim(W^s(g\cdot x)\cap
  M')+\dim(W^u(x)\cap M')\leq \dim M'$ any transversal intersection
  between $W^u(x)$ and $W^s(g\cdot x)$ must be trivial, i.e.\ the
  connection cannot be robust; this proves (i).
  
  To see (ii), note that $g\cdot x\in \Fix(K)$ implies that $x\in
  g^{-1}\cdot \Fix(K)$ and so $x\in \Fix(K^{g^{-1}})$. By (i) $K\neq
  K^{g^{-1}}$ and so the isotropy $H$ of $x$ contains both $K$ and
  $K^{g^{-1}}$ but is not equal to either.
\end{proof}

\begin{corollary}
  If $G$ is abelian there can be no robust \RHT s.
\end{corollary}

\begin{proof}
  This is because if $G$ is abelian then $N_G(K)=G$ for all $K<G$ and
  so there is no twist satisfying Proposition~\ref{proptwist} (i).
\end{proof}

In Proposition \ref{prop:D_n} we show similarly that there are no
robust \rht s for actions of the dihedral group $\D_n$, unless $n$ is
a multiple of 4.

\begin{remark}
  Although the previous corollary excludes the possibility of robust
  \RHT s for Abelian groups, one can find robust relative
  \emph{heteroclinic} cycles for Abelian groups. For example, one can
  break the cyclic symmetry of the Guckenheimer-Holmes cycle by making
  $\lambda$ in (\ref{eqGHmodel}) vary with index. This breaks the
  symmetry to the Abelian group $(\Z_2)^3$ while leaving the same
  cycle robust. However, such perturbations will break the relative
  equilibria into three families of relative equilibria and the
  connections will no longer be homoclinic. Work of Melbourne {\em et
    al.} in particular gives general methods for locating robust
  relative heteroclinic cycles in terms of cycles in the isotropy
  lattice and conditions on the isotypic decompositions of certain
  isotropy subgroups \cite{MCG89}.
\end{remark}

%%%%%%%%%%%%%%%%%%%%%%%%%%%%%%%%

\section{Classification of homoclinic triples}
\label{sec:classification} \setcounter{equation}{0}

Motivated by Proposition \ref{proptwist} we make the following
definition.  Given a group $G$, a \emph{homoclinic triple} for $G$ is
a triple $(K,g,H)$, where $K$ and $H$ are subgroups of $G$ and $g\in
G$, satisfying
\begin{description}
\item[HT1:] $gH\cap N_G(K)=\emptyset$ and
\item[HT2:] $K < H\cap H^g$.
\end{description}

Note that [HT2] is equivalent to assuming that $H>K\cup K^{g^{-1}}$.
Homoclinic triples in $G$ do not always give rise
to robust \rht s of $G$-equivariant vector fields, 
though we show in Section \ref{secconst} that with an
additional hypothesis they do.

It is clear that if $(K,g,H)$ is a homoclinic triple for $G$, and
$G<G'$ then $(K,g,H)$ is a homoclinic triple for $G'$.

Two homoclinic triples $(K,g,H)$ and $(K',g',H')$ are said to be
\emph{conjugate} if there is an element $f\in G$ such that $K'=K^f$,
$H'=H^f$ and $g'=g^f=fgf^{-1}$.  Moreover if there is a vector field
such that the the first is the triple associated to an \rht\ $\gamma$
then the second is the triple associated to the \rht\ $f\cdot\gamma$,
for the same vector field.

A further equivalence can be given by the relation $(K,g,H) \sim
(K,g^{-1}, H^g)$.  Indeed, it is easy to see that if a vector field
$v$ gives a robust \RHT\ with triple $(K,g,H)$, then the opposite
vector field $-v$ gives a robust \RHT\ with triple $(K,g^{-1}, H^g)$.

What seems less obvious is the following fact, whose relationship with
the dynamics is not clear.

\begin{proposition} \label{prop:nu}
  If $(K,g,H)$ is a homoclinic triple, then so is $(K,\nu g,H)$ for
  all $\nu\in N_G(K)$.
\end{proposition}

\begin{proof}
  Suppose $(K,g,H)$ is a homoclinic triple, and let $g'=\nu g$ for
  $\nu\in N_G(K)$. Then firstly if $gH\cap N_G(K)=\emptyset$ then $\nu
  gH\cap \nu N_G(K)=\emptyset$, but $\nu N_G(K) = N_G(K)$, and
  secondly we have $H>K$ and since $\Kginv = K^{(\nu g)^{-1}}$ we
  have $H>K^{g'^{-1}}$ as required. 
\end{proof}

It is clear that if $(K,g,H)$ is a homoclinic triple for $G$, and
$G_1<G$ is a subgroup containing $g$ then $(K\cap G_1, g, H\cap G_1)$
is a homoclinic triple for $G_1$.  More generally, the property of
being a homoclinic triple is preserved under pull-back:

\begin{proposition} \label{prop:HM}
  Let $\phi: G_1\to G_2$ be a homomorphism, and let $(K_2,g_2,H_2)$ be
  a homoclinic triple in $G_2$, with $g_2\in \mathrm{Im}(\phi)$.
  Define $K_1=\phi^{-1}(K_2)$, $H_1=\phi^{-1}(H_2)$ and let
  $g_1\in\phi^{-1}(g_2)$.  Then $(K_1,g_1,H_1)$ is a homoclinic triple
  in $G_1$.
\end{proposition}

\begin{proof}
  It is easy to check that if $g_1H_1\cap N_{G_1}(K_1)\neq\emptyset$
  then the image under $\phi$ of any element in this intersection
  belongs to $g_2H_2\cap N_{G_2}(K_2)$ which by hypothesis is empty,
  thus establishing [HT1].  Similarly, [HT2] also follows simply by
applying $\phi^{-1}$.
\end{proof}

A useful result for classifying homoclinic triples is the following:

\begin{proposition}
  Suppose $(K,g,H)$ is a homoclinic triple.  Then so is $(K,g,H_0)$,
  where $H_0$ is the group generated by $K$ and $\Kginv$.
\end{proposition}

\begin{proof}
  It is clear that [HT1] is still satisfied while [HT2] is satisfied
  by construction.
\end{proof}

Such triples $(K,g,H_0)$ are called \emph{minimal} homoclinic triples,
since necessarily $H_0<H$.

The aim is to classify such minimal homoclinic triples up to the
equivalence relation generated by conjugacy and the equivalence of
Proposition \ref{prop:nu}.  The procedure is as follows.  Fix a
subgroup $K$, and choose a distinct conjugate $K'=\Kginv$, for some
$g$ which is well-defined modulo $N_G(K)$ (as in the proposition
above).  Let $H=H_0$ be the group generated by $K$ and $K'$ (the
smallest possible group satisfying [HT2$'$]), and finally check the
remaining criterion [HT1].  In simple cases---such as small
permutation groups---the checking can be done by \textsc{Maple}.

Having found the smallest possible $H_0$ for a given pair $(K,g)$ one
can then consider enlarging $H$ until [HT1] is no longer satisfied.

The smallest non-abelian group is the symmetric group $S_3=\D_3$, but there
is only one non-normal subgroup of $S_3$ up to conjugation, namely a
$\Z_2$, and for this subgroup there is no possible twist, as is
readily checked.  

One of the next smallest non-abelian group is the symmetry
group $\D_4$ of the square, and for this there are homoclinic triples,
as found in the Kuramoto-Sivashinsky example. In fact there are two
homoclinic triples, which are not conjugate, though they are
equivalent under an \emph{outer} automorphism of the group.

Another non-abelian example of the same size is the eight-element
quaternion group $Q=\{\pm 1,\pm i, \pm j, \pm k\}$ where
$i^2=j^2=k^2=ijk=-1$. The only nontrivial groups $K$ and $H$ such that
$G>H>K$ are nontrivial containments have $K=\{\pm 1\}$. As $N_G(K)=G$
we cannot satisfy HT1 and so $Q$ supports no homoclinic triples.

In each of the tables of the sections that follow, the final column
states whether $K$ is a normal subgroup of both $H$ and $H^g$, or
equivalently whether $H$ is contained in the normalizers of both $K$
and $\Kginv$.  This is a condition we will use in the existence
theorem of Section \ref{secconst}.

%%%%%%%%%%%%%%%%%%%%%%%%%%%
\subsection{Dihedral groups}
\label{sec:dihedral}

Denoted $\D_n$, these are the symmetry groups of the regular polygons.
In order to determine homoclinic triples in $\D_n$ we introduce some
notation, and recall a few basic properties of $\D_n$.  Let $\rho$ be
the rotation by $2\pi/n$, and $\kappa$ a reflection, so that
$\D_n = \left<\kappa,\rho\right>$, and $\rho\kappa = \kappa\rho^{-1}$.
The elements of $\D_n$ can then be written as
$$
\{1,\rho,\rho^2,\dots,\rho^{n-1},\; \kappa, \kappa\rho, \dots,
\kappa\rho^{n-1}\}.
$$
Any subgroup of the cyclic subgroup $\Z_n$ of
$\D_n$ is normal in $\D_n$, so candidates for the symmetry group $K$
of a robust \rht\ must contain a reflection.  If $n$ is even then there
are two distinct conjugacy classes of reflection, given by $\kappa$ and
$\kappa\rho$.  More generally, $\kappa\rho^r \sim \kappa\rho^s$ iff
$r$ and $s$ are of the same parity modulo $n$.

%%%%%%%%%%%%%%%%%%%%%%%%%%%
\paragraph{$G=\D_4$}
\begin{center}
  \begin{tabular}{l|l|l|l|l|c}
    $K$ \qquad& generators of $K$ \qquad& twist, $g$ \quad& 
    $H$ & generators of $H$  & $K\lhd H,H^g$\cr
    \hline
    $\Z_2$ & $\kappa$ & $\rho$ & $\D_2$ & 
    $\kappa,\rho^2$ &  Yes \cr
    
    $\Z_2$ & $\kappa\rho$ & $\rho$ & $\D_2'$ & 
    $\kappa\rho,\rho^2$ & Yes\cr
    
    \hline     
  \end{tabular}
\end{center}
  
One of these two homoclinic triples, say with $K$ generated by
$\kappa$, occurs in the Kuramoto-Sivashinsky system.  The other one
does not as $\D_2'$ is not an isotropy subgroup.  It should be pointed
out that the two rows are equivalent under an \emph{outer}
automorphism of the group but they are not conjugate, which is why
both are included.
%%%%%%%%%%%%%%%%%%%%%%%

\newpage

%%%%%%%%%%%%%%%%%%%%%%%
\paragraph{$G=\D_8$}
This group contains $\D_4$ so two cases follow from those of $\D_4$,
namely those with twist $\rho^2$.
\begin{center}
  \begin{tabular}{l|l|l|l|l|c}
    $K$ \qquad& generators of $K$ \qquad& twist, $g$ \quad& 
    $H$ & generators of $H$ &  $K\lhd H,H^g$ \cr    
    \hline
    $\Z_2$ & $\kappa$ & $\rho$ & $\D_4$ & 
    $\kappa,\rho^2$ &  No \cr
    
    $\Z_2$ & $\kappa\rho$ & $\rho$ & $\D_4$ & 
    $\kappa\rho,\rho^2$ &  No \cr
    
    $\Z_2$ & $\kappa$ & $\rho^2$ & $\D_2$ & 
    $\kappa,\rho^4$ & Yes \cr
    
    $\Z_2$ & $\kappa\rho$ & $\rho^2$ & $\D_2$ & 
    $\kappa\rho,\rho^4$ & Yes \cr
    
    $\D_2$ & $\kappa,\rho^4$ &  $\rho$ & $\D_4$ &
    $\kappa,\rho^2$ &  Yes \cr
    
    $\D_2$ & $\kappa\rho,\rho^4$ &  $\rho$ & $\D_4$ &
    $\kappa\rho,\rho^2$  &  Yes \cr
    
    \hline
  \end{tabular}
\end{center}
The first two lines are equivalent under an \emph{outer} automorphism
of $G$, as are each of the next pair and the final pair.

\begin{proposition} \label{prop:D_n}
  $\D_n$ contains homoclinic triples iff $n$ is a multiple of 4.
\end{proposition}

In particular, by Proposition \ref{proptwist} there are no robust \rht
s in systems with dihedral symmetry $\D_n$ unless $n$ is a multiple of
4.

\begin{proof}
  First, if $n$ is a multiple of 4, then $\D_4 < \D_n$, so the result
  follows from finding a homoclinic triple for $\D_4$. This is given
  in the table above.
  
  For the converse, let $K<G$ be a non-normal subgroup, and suppose
  $(K,g,H)$ is a homoclinic triple. As argued above, we may assume
  $\kappa \in K$.  Without loss of generality we can take the twist to
  be $g=\rho^r$ for some $r$, for by Proposition \ref{prop:nu} $(K,
  \rho^r, H)$ is a homoclinic triple if and only if $(K, \kappa\rho^r,
  H)$ is since $\kappa\in K < N_G(K)$.  Then, $\Kginv$ contains
  $\rho^{-r}\kappa\rho^r = \kappa\rho^{2r}$.  It follows in
  particular that $h=\rho^{2r}\in H$.
  
  Suppose first that $n$ is odd, and put $a=(n-1)/2$.  Then $gh^a =
  \rho^{r(2a+1)}= \rho^{rn}=1\in N_G(K)$ contradicting the hypothesis
  that $(K,g,H)$ is a homoclinic triple.
  
  Now suppose $n=2p$, with $p$ odd.  Then $\rho^p$ (rotation by $\pi$
  in $\D_n$) is in the centre of $\D_n$ so that it belongs to $N_G(K)$
  for any subgroup $K$. Now repeat the argument above but with
  $a=(p-1)/2$, and one finds that $gh^a =\rho^p\in N_G(K)$.
\end{proof}

We now give a complete classification of homoclinic triples for
subgroups of $\D_{4n}$.  For uniformity of notation, we write
$\Z_2=\D_1$.

\begin{theorem} \label{thm:D4n}
  Let $G=\D_{4n}$.  Up to equivalence, the only homoclinic triples are
  of the form $(K,g,H)=(\D_q, \rho^r, \D_{p})$, satisfying,
  \begin{description} 
     \item[(i)] $q\vert n$ 
     \item[(ii)] $t(r)\leq t(n/q)$ 
     \item[(iii)] $q\vert p$, $q\neq p$, $p\vert (2n)$ and $r\equiv
       (2n/p) \mathrm{\ mod\ }(4n/p)$
  \end{description}
  Here $t(r)$ is the multiplicity of 2 in the prime decomposition of
  an integer $r$.  Moreover all such triples are homoclinic triples.
\end{theorem}

\begin{proof}
As has already been pointed out, homoclinic triples cannot have
$K=\C_q$ for these are normal subgroups of $G=\D_n$. Thus we can
suppose $K=\D_q$ for some $q\geq1$ which necessarily divides
$4n$. Write $s=4n/q$, and let $\kappa$ be a reflection in $K$. Thus
$\D_q=\left<\kappa,\rho^s\right>$.  The normalizer $N_G(K)$ is
$\D_{q'}$, where $q'=q$ if $s$ is odd, and $q'=2q$ if $s$ is
even. Thus, $N_G(K) = \left<\kappa,\rho^{s'}\right>$, where $q's'=4n$.
If $K=\left<\kappa\right>$ (ie, if $q=1$) then $s=4n$ and $s'=2n$.

Without loss of generality, we can assume the twist $g=\rho^r$ for
some $r$, for the alternative form is $\kappa\rho^r$, but the two are
equivalent by Proposition \ref{prop:nu}.  Furthermore, by the same
equivalence, we can choose $0 < r< s'$.  

Let $H_0$ be the subgroup generated by $K$ and $\Kginv$. Then
$$ H_0 = \left<\kappa,\,\rho^{2r},\,\rho^s\right>,$$
and
$$ gH_0 = \{\rho^{r+k2r+\ell s},
      \kappa\rho^{-r+k2r+\ell s} \mid k,\ell\in\Z \},$$ 
so that $gH_0\cap N_G(K)\neq\emptyset$ if and only if there are
integers $j,k,\ell$ such that 
$$\rho^{r+js+2kr} =\rho^{\ell s'} \quad \mbox{or} \quad
\kappa\rho^{-r+js+2kr} = \kappa\rho^{\ell s'}.$$
This is equivalent to there being integers $k,\ell$ such that
\begin{equation} \label{eq:D4n-cond}
(2k+1)r  = \ell s'
\end{equation}
If $q$ does not divide $n$ then $s'$ is odd and there are always
solutions to this equation (for all $r$), so establishing (i).

Now suppose $q\vert n$, so that $s'$ is even.  Then
(\ref{eq:D4n-cond}) has solutions iff $t(r)\geq t(s')$ so establishing
(ii).

Finally, let $H=\D_p = \left<\kappa,\,\rho^a\right>$, for some $a$
which must divide both $s=2s'$ and $2r$ (so that $H$ contains both $K$
and $\Kginv$). If $a$ divides $r$ then $g\in H$ which is not
possible. Thus $a$ is even since it divides $2r$, and $r\equiv (a/2)
\textrm{\ mod\ }a$.  

In this case $gH \cap N_G(K)\neq\emptyset$ if and only if there are
integers $k,\ell$ such that 
$$r + ka = \ell s'.$$ %
Since $a\vert s$ and $s=2s'$ we have either $s'\equiv 0 \textrm{\ mod\
}a$ or $s'\equiv (a/2) \textrm{\ mod\ }a$.  In the first case there
are no solutions to the equation, while in the second there are
solutions, so establishing (iii).  \end{proof}

If $K=\D_q$ then, up to equivalence, the only case where $K$ and
$\Kginv$ are normal subgroups of $H$ is $H=\D_{2q}$ and
$g=\rho^{n/q}$.

%%%%%%%%%%%%%%%%%%%%%%%%%%%
\subsection{Symmetric groups}
\label{sec:S_n}

\paragraph{$G=S_3\simeq \D_3$} 
We have already proved that this group has no homoclinic triples
(Proposition \ref{prop:D_n}).

%%%%%%%%%%%%%%%%%%%%%%%%%%%
\paragraph{$G=S_4\simeq\mathbb{T}_d$} 
The only inequivalant homoclinic triples are given in the following
table:

\begin{center}
  
  \begin{tabular}{l|l|l|l|l|c}
    $K$ \qquad& generators of $K$ \qquad& twist, $g$ \quad& 
    $H$ & generators of $H$ &  $K\lhd H,H^g$ \cr    
    \hline
    $\Z_2$   & $(1\, 2)$          & $(1\, 3)(2\, 4)$ & 
    $\D_2$ & $(1\,2), (3\,4)$  & Yes\cr
    $\Z_2$   & $(1\, 2)(3\, 4)$   & $(1\,3)$ & 
    $V_4$ & $(1\, 2)(3\, 4),(1\, 4)(2\, 3)$ & Yes\cr
    \hline
  \end{tabular}
\end{center}

%%%%%%%%%%%%%%%%%%%%%%%%%%%
\paragraph{$G=S_5$} 
the only homoclinic triples appear to be those that come from the
inclusion $S_4<S_5$.

%%%%%%%%%%%%%%%%%%%%%%%%%%%
\paragraph{$G=S_6$}
 This list is by no means expected to be exhaustive; there are many
different classes of subgroup of $S_6$. 
\begin{center}
\begin{tabular}{l|l|l|l|l|c}
    $K$ \qquad& generators of $K$ \qquad& twist, $g$ \quad& 
    $H$ & generators of $H$ &  $K\lhd H,H^g$ \cr    
    \hline
    $\Z_2$   & $(1\, 2)$          & $(1\, 3)(2\, 4)$ & 
      $\D_2$ & $(1\,2), (3\,4)$ & Yes\cr
    
    $\Z_2$   & $(1\, 2)(3\,4)$    & $(1\, 3)$ & 
      $V_4$ & $(1\, 2)(3\,4),(1\, 3)(2\,4) $ & Yes \cr

    $\Z_2$   & $(1\, 2)(3\,4)$    & $(1\, 5)(2\, 6\, 3)$ & 
      $\D_4$ & $(1\, 2)(3\,4),(2\, 4)(5\, 6) $ & No \cr

    $\Z_2^2$   & $(1\, 2),(3\,4)$   & $(1\, 5)(2\, 6) $ &
      $\Z_2^3$ & $(1\, 2),(3\,4), (5\, 6) $ & Yes  \cr
      
    $\Z_2^2$   & $(1\, 2),(3\,4)$   & $(1\, 3\, 5)(2\, 6) $ &
      $\Z_2\times S_4$ & $(1\, 2),(3\,4), (3\,5),(3\,6) $ & No  \cr
     
    $\C_3$     & $(1\, 2\, 3)$ &    $(1\, 4)(2 \,5)(3\, 6)$
      &$\Z_3^2$ &  $(1\, 2\, 3),(4\, 5\, 6)$ & Yes \cr

    $\C_6$   & $(1\, 2\, 3\, 4\, 5\, 6)$ &    $(2\, 4)$
      & W &  $(1\, 2\, 3\, 4\, 5\, 6), (1\, 2\, 3\, 6\, 5\, 4)$ & No \cr

    $\D_6$   & $(1\, 6)(2\, 5)(3\, 4), (1\, 2\, 3\, 4\, 5\, 6)$ &  $(2\, 4)$
      & W &  $(1\, 2\, 3\, 4\, 5\, 6), (1\, 2\, 3\, 6\, 5\, 4)$ & No \cr

    $\Z_6$     & $(1\, 2\, 3),(4\, 5)$ &    $(1\, 4)(2 \,5\, 3\, 6)$
      &$S_3^2$ &  $(1\, 2),(1\, 3),(4\, 5),(4\, 6)$ & No \cr

    $S_3$     & $(1\, 2),(1\, 3)$ &    $(1\, 4)(2 \,5)(3\, 6)$
      &$S_3^2$ &  $(1\, 2),(1\, 3),(4\, 5),(4\, 6)$ & Yes \cr

    \hline
  \end{tabular}
The group $W$ has order 36 and has the two generators as shown. 
\end{center}

%%%%%%%%%%%%%%%%%%%%%%%%%%%
\subsection{Wreath products} 
\label{sec:wreath}

Wreath products are the natural form of symmetry group occurring in
systems of coupled cells. They are groups of the form $G =
\mathcal{L}\wr\mathcal{G}$ where $\mathcal{G}<S_n$ is a subgroup of a
permutation group and $\mathcal{L}$ is a nontrivial compact group; see
for example \cite{Dio&al96,Dia&al99}.  Here we consider two specific
examples.

%%%%%%%%%%%%%%%%%%%%%%%%%%%
\paragraph{$G=\Z_2\wr\Z_3\simeq A_4\times \Z_2\simeq\mathbb{T}_h$} 
This is the group occurring in the Guckenheimer-Holmes example, see
\S\ref{secGHmodel}.  We write $\rho$ for the generator of $\Z_3$ and
$\kappa_j$ for the generator of the $j^{\rm th}$ copy of $\Z_2$.

\begin{center}
  \begin{tabular}{l|l|l|l|l|c}
    $K$ \qquad& generators of $K$ \qquad& twist, $g$ \quad& 
    $H$ & generators of $H$ & $K\lhd H,H^g$ \cr
    \hline
    $\Z_2$ & $\kappa_1$ & $\rho$ & $\D_2$ & 
    $\kappa_1,\kappa_2$ & Yes \cr
    
    $\Z_2$ & $\kappa_1\kappa_2$ & $\rho$ & $\D_2$ & 
    $\kappa_1\kappa_2, \kappa_1\kappa_3$ & Yes \cr
    
    $\D_2$ & $\kappa_1,\kappa_2$ &  $\rho$ & $\Z_2^3$ &
    $\kappa_1,\kappa_2,\kappa_3$ & Yes \cr
    
    $\Z_2^2$ & $\kappa_1,\kappa_2\kappa_3$ & $\rho$ & $\Z_2^3$ &
    $\kappa_1,\kappa_2,\kappa_3$ & Yes\cr
    \hline     
  \end{tabular}
\end{center}

The first row is the case occurring in the Guckenheimer-Holmes
example. Observe that this is only one of a number of possible
homoclinic triples.

%\newpage

%%%%%%%%%%%%%%%%%%%%%%%%%%%

\paragraph{$G=\Z_2\wr S_3$} This is the octahedral group $\mathbb{O}_h$.
Write $\rho=(1\,2\,3)$ and $\sigma=(1\,2)$ in $S_3$ and let $\kappa_j$
be the generator of $\Z_2$ in the $j^{\mbox{th}}$ position.

\begin{center}
  \begin{tabular}{l|l|l|l|l|cr}
    $K$ \quad\ & generators of $K$ & twist $g$ & $H$ & generators of $H$ & $K\lhd H,H^g$ \cr
    \hline
    $\Z_2$ & $\kappa_1$ & $\sigma$ & $\Z_2^2$ & 
    $\kappa_1, \kappa_2$ & Yes \cr
        
    $\Z_2$ & $\kappa_1\kappa_2$ & $\sigma$ & $\Z_2^2$ & 
    $\kappa_1\kappa_2, \kappa_1\kappa_3$ & Yes \cr

    $\Z_2$ & $\sigma$ & $\kappa_1$ & $\D_2$ & 
    $\sigma, \kappa_1\kappa_2$ & Yes \cr
    
    $\Z_2$ & $\sigma\kappa_3$ & $\kappa_1$ & $\D_2$ & 
    $\sigma\kappa_3, \kappa_1\kappa_2\sigma$ & Yes \cr

    $\Z_2^2$ & $\sigma,\kappa_3$ & $\kappa_1$ & $ \Z_2^3$ & 
    $\sigma,\kappa_1\kappa_2, \kappa_3$ & Yes \cr

    $\Z_2^2$ & $\kappa_1,\kappa_2$ & $\rho$ & $\Z_2^3$ & 
    $\kappa_1, \kappa_2, \kappa_3$ & Yes \cr
    \hline
  \end{tabular}
\end{center}

%%%%%%%%%%%%%%%%%%%%%%%%%%%

\begin{example} \label{eg:wreath}
  It is easy to extend some of these to the general case
  $G=\mathcal{L}\wr \mathcal{G}$.  We set $K$ to be a subgroup
  consisting of $\mathcal{L}$ in only one component, say
  $\mathcal{L}_1$ and identity elsewhere.  If we take
  $g\in\mathcal{G}$ which takes cell 1 to say, cell 2 (eg $g=(1\,
  2)\in S_n$ corresponding to $\sigma$ in the table above or $g=(1\,
  2\, 3)$ corresponding to $\rho$ in the previous table) then
  $N_G(K)=\mathcal{L}_1 \times (\mathcal{L}\wr\mathcal{G'})$ for some
  subgroup $\mathcal{G'}$ of $\mathcal{G}$, and $K^{g^{-1}} =
  \mathcal{L}_2$.  Let $H=\mathcal{L}_1\times \mathcal{L}_2$
  (or any other subgroup of $N_G(K)$ containing $\mathcal{L}_1\times
  \mathcal{L}_2$) then in addition to being a homoclinic triple, one
  also finds that $K\lhd H,H^g$.
\end{example}

%%%%%%%%%%%%%
\section{Construction of robust {\sc rht}s}
\label{secconst} \setcounter{equation}{0}

Proposition~\ref{proptwist} gives group-theoretic condition on the
twist of an \RHT\ necessary for it to be robust.  We now give
sufficient conditions that allow construction of robust \RHT s.

\begin{theorem}\label{thmsufficient}
  Let $G$ be a compact Lie group acting on $M$.  Suppose $K$ and $H$
  are isotropy subgroups of this action and $(K,g,H)$ is a homoclinic
  triple satisfying the two further conditions (one `local' the other
  `global'):
  \begin{description}
  \item[HTL] $H < N_G(K) \cap N_G(K^{g^{-1}})$,
  \item[HTG] there is a point $x\in M$ with isotropy $H$ and a
    continuous path of points with isotropy $K$ that joins $x$ to
    $g\cdot x \in \Fix(H^g)$.
  \end{description}
  Then there exists a non-empty open set of equivariant vector fields
  on $M$ with a robust \RHT\ of isotropy $K$, based on a point in
  $\Fix(H)$ and with twist $g$.
\end{theorem}

\begin{proof}
  Note that if $K$ is an isotropy subgroup then so is $K^g$ and by the
  fact that $(K,g,H)$ is a homoclinic triple, $K$ and $\Kginv$ are
  distinct isotropy subgroups and so have distinct fixed point
  subspaces.  Let $X_1$ be the open set of equivariant vector fields
  on $M$ such that there is a normally hyperbolic (relative)
  equilibrium point $x$ with isotropy $H$.
  
  Consider the action of $H$ on $T_xM$, and consider the three
  subspaces $T_H=\Fix(H,T_xM)$, $T_K=\Fix(K,T_xM)$ and
  $T_K'=\Fix(K^{g^{-1}},T_xM)$. By [HTL] these are each $H$-invariant
  subspaces.  It follows that there is an $H$-invariant decomposition
  of $T_xM$ as
  $$
  T_xM = T_H\oplus T_0\oplus T_1\oplus T_2 \oplus W
  $$
  where
  \begin{eqnarray*}
    T_H\oplus T_0 &=& T_K \cap T_K' \\
    T_H\oplus T_0 \oplus T_1 &=& T_K\\
    T_H\oplus T_0 \oplus T_2 &=& T_K'.
  \end{eqnarray*}
  The fact that $T_H\subset T_K\cap T_K'$ follows from the homoclinic
  triple property.  $T_1$ and $T_2$ are of course of the same
  dimension.  We treat explicitly the case where $G$ is finite; the
  general case can be deduced by intersecting the above decomposition
  with a local slice to the group orbit through $x$.
  
  While the linearization $u\mapsto Lu$ of the vector field at $x$
  leaves invariant each of $T_H, T_K, T_K'$ and $T_K\cap T_K'$, it
  does not necessarily leave invariant the entire decomposition.
  However, there are eigenvalues associated to each subspace.  For
  example since $T_H$ and $T_H\oplus T_0$ are invariant, it follows
  that the linear vector field descends to a vector field on
  $(T_H\oplus T_0)/T_H \simeq T_0$ whose eigenvalues we call
  $\sigma_0(L)$.  Similar constructions define $\sigma_1(L),
  \sigma_2(L)$ etc.  
  
  Define the class $X_2$ of vector fields in $X_1$ such that
  $$
  m_1^+ + m_2^- > \dim T_1,
  $$
  where $m_1^+$ is the number of eigenvalues in $\sigma_1(L)$ with
  positive real part, counting multiplicity, and similarly $m_2^-$
  counts the number of eigenvalues in $\sigma_2(L)$ with
  \emph{negative} real part.  This class $X_2$ is clearly open in
  $X_1$.  Furthermore, the members of this class satisfy
  \begin{equation} \label{eq:dimensions}
    (W^u(x)\cap \Fix(K))+\dim(W^s(g\cdot x)\cap\Fix(K))>\dim \Fix(K).
  \end{equation}
  This is because $g\cdot\left(W^s(x)\cap \Fix(K^{g^{-1}})\right) =
  W^s(g\cdot x)\cap\Fix(K)$.

  Now consider any continuous path $\gamma(t)$ with finite arc length,
  connecting $x$ to $g\cdot x$ in $\Fix(K)$ such that $\gamma(t)$ has
  isotropy $K$.  Given a tubular neighbourhood $N$ of $\gamma$ such
  that $g(N)$ does not intersect $N$, there is an open set $X_3\subset
  X_2$ of vector fields such that $W^u(x)$ and $W^s(g\cdot x)$ intersect on
  a path within $N$.
  
  Since $\dim(W^u(x)\cap \Fix(K))+\dim(W^s(x)\cap\Fix(K))\geq \dim
  \Fix(K)+1$ there is an open set $X_4\subset X_3$ of flows that have
  transverse intersection of these manifolds and hence the \RHT s in the
  open set $X_4$ are robust as required.
\end{proof}

\begin{remark}
  Notice that although $K\neq K^g$ we do not require that $H$ and
  $H^g$ differ for this result. The examples in
  Sections~\ref{secGHmodel} and \ref{secKSmodel} have $H\neq H^g$ and
  $H=H^g$ respectively.  Moreover the constructed \RHT s may or may not
  be attractors.
\end{remark}

\begin{remark} \label{rmk:HTL}
  Condition [HTL] is not necessary, but it simplifies the proof
  considerably. Any more general hypothesis on the representation of
  $H$ on $T_xM$ would suffice provided it allows construction of open
  sets of vector fields satisfying equation (\ref{eq:dimensions}); see
  \S\ref{sec:D_8-example} for an example. It is worth noticing that
  [HTL] is satisfied whenever $H$ is abelian, which holds in almost
  all the examples in Section \ref{sec:classification}.  
\end{remark}

\begin{remark}
  Condition [HTG] is obviously necessary, and requires
  knowledge of the group action in order to be verified.  However, if
  the action is a complex or symplectic representation of a finite
  group then the symmetry-type sets (ie those with constant isotropy)
  are always connected, so that [HTG] is always verified; see Theorem
  \ref{thmsufffinite}.  
\end{remark}

%%%%%%%%%%%

\subsection{An example of a robust {\sc rht}  with $SO(3)$ symmetry.}

Theorem~\ref{thmsufficient} can be applied to produce new examples of
robust \RHT s in a variety of contexts. For example, let $G=SO(3)$ and
let $V_1$ be the 3-dimensional irreducible representation (irrep) of
$SO(3)$ ($V_1=\R^3$), and $V_2$ the 5-dimensional irrep, consisting of
trace-0 $3\times3$ real symmetric matrices (equivalently, the first
and second order spherical harmonics respectively). Let $M =
(V_1\oplus V_2)\otimes\C$. This is a 16-dimensional representation of
$SO(3)$, isomorphic to the sum $\C^3\oplus V_2^\C$, where $V_2^\C$
consists of \emph{complex} trace-0 symmetric $3\times3$ matrices.
$SO(3)$ acts on $V_1$ by multiplication on the left, and on $V_2$ by
conjugation.

Note that action of $SO(3)$ on $V_2$ has non-trivial generic orbit
type (equal to the conjugacy class of the subgroup $H$ defined below),
since every symmetric matrix is diagonalisable.  On the other hand, if
the real and imaginary parts of a complex matrix in $V_2^\C$ have no
common eigenvector then the isotropy of the point is trivial.

Denote by $R^x_\theta$ the rotation by $\theta$ about the $x$-axis.
Let $H=\D_2$ be the group generated by $R^x_\pi$ and $R^y_\pi$.  Then
$\Fix(H)$ consists of diagonal matrices.  Let $K =
\left<R^x_\pi\right>$, so that $\Fix(K)$ consists of certain block
diagonal matrices, and let $g=R^y_{\pi/2}$.  

\begin{theorem} 
  Suppose $G=SO(3)$ and $M$ is the 16-dimensional representation as
  above. Then there exist robust \RHT s on $M$ for vector fields with
  symmetry $G$.
\end{theorem}

\begin{proof}
  We verify that the hypotheses of Theorem~\ref{thmsufficient} are
  satisfied and then use the conclusion of that theorem. To see that
  $(K,g,H)$ is a homoclinic triple, note that $g$ does not fix
  $\Fix(K)$ as it maps the $x$ axis onto the $z$ axis, and observe
  that $H$ is generated by $K$ and $K^{g^{-1}}$.  For [HTL] it
  suffices to point out that $N_G(K)=N_G(K^{g^{-1}})=O(2)$ (generated
  by all $R^x_\theta$ and $R^y_\pi$), and this contains $H$.
  
  There remains to show that the global connectivity condition is
  satisfied, which can be done explicitly.  Since $SO(3)$ acts
  independently on each summand in $M$, it follows that for any
  subgroup
  $$\Fix(K,M) = \Fix(K,V_1^\C) \oplus \Fix(K,V_2^\C),$$ %
  and moreover that $\Fix(K,V^\C) = \Fix(K,V)^\C$, and finally that
  the isotropy subgroup of $(u+iv,A+iB)\in M$ is the intersection of the
  isotropy subgroups of each of $u,v,A$ and $B$.
  
  Thus $\Fix(K,M)$ consists of points of the form
  \begin{equation}
    \label{eq:fixK}
    \left[\left(
        \begin{array}{c} 
          t\\ \hline 0     
        \end{array}\right),\;
      \left(
        \begin{array}{c|c}
          a&0\\ \hline 0 &A
        \end{array}
      \right)\right]
  \end{equation}
  where $a,t\in\C$ and $A$ is a $2\times2$ complex symmetric matrix of
  trace $\mathop{\rm tr}(A)=-a$.  A typical point in $\Fix(H,M)$ is a
  point of the form $U=(0,\mathop{\rm diag}[\alpha,\beta,\gamma])$,
  and then $g\cdot U = (0,\mathop{\rm diag}[\gamma,\beta,\alpha])$. We
  need to show that there is a path in $\Fix(K,M)$ connecting $U$ and
  $g\cdot U$ consisting entirely of points with isotropy precisely
  $K$.
  
  Now, $\dim\Fix(K,M)=8$, and $\dim\Fix(H)=4$. Furthermore $N_G(K)
  \simeq O(2)<SO(3)$, so that the set of points in $\Fix(K,M)$ with
  orbit type $(H)$ is of dimension $4+1=5$ and so its complement in
  $\Fix(K)$ is connected.  The only other points in $\Fix(K,M)$ with
  higher symmetry are those fixed by $SO(2)<N_G(K)$ (generated by
  rotations about the $x$-axis) which are just those of the form
  (\ref{eq:fixK}) with $A$ a multiple of the identity, which is a set
  of real dimension only 4, and so again does not separate
  $\Fix(K,M)$.
\end{proof}

%%%%%%%%%%%%%
\subsection{A class of finite groups admitting robust {\sc rht}s}

The next result gives sufficient conditions on a finite group $G$ such
that a representation $M$ of $G$ can be found that admits robust 
\RHT s for $G$ equivariant vector fields on $M$.

\begin{theorem}\label{thmsufffinite}
  Suppose that $G$ is a finite group with a homoclinic triple
  $(K,g,H)$ satisfying in addition [HTL].  Let $M$ be the complex
  regular representation of $G$ (of real dimension $2|G|$).  Then
  there is an open set of $G$-equivariant vector fields on $M$ for
  which there is an \RHT\ with isotropy $K$ and twist $g$.
\end{theorem}

\begin{proof}
  This follows from Theorem~\ref{thmsufficient}. Indeed, every
  subgroup of $G$ is an isotropy subgroup for this representation so
  that $K$ and $K^{g^{-1}}$ and $H$ are isotropy subgroups as
  required, and hypothesis [HTL] is satisfied by assumption.  Since
  all strata are even-dimensional there are no connectivity
  restrictions and we can satisfy [HTG].    
\end{proof}

\begin{example} 
  As an example of an application of this theorem, consider the
  eight-element dihedral group
  $$
  \D_4=\langle \{
  \rho,\kappa~:~\rho^4=\kappa^2=1,~~\rho\kappa=\kappa\rho^{-1}
  \}\rangle.
  $$
  This has a subgroup $K=\langle \kappa \rangle$ and element
  $g=\rho$ that satisfy the hypotheses of the theorem (in this case
  $H=\langle \kappa, \rho^2\rangle$) and so we can construct robust
  cycles with $\D_4$ symmetry for a flow on $\R^{16}$. This is
  identical in form to the Kuramoto-Sivashinsky example in
  Section~\ref{secKSmodel} but in a much higher dimensional space.
\end{example}

\begin{example}
  Similarly, one can use Theorem~\ref{thmsufffinite} to show that the
  group occurring in the Guckenheimer-Holmes example in
  Section~\ref{secGHmodel} permits robust \RHT s. Namely, let
  $G=\Z_2\wr S_3$, see \S\ref{sec:wreath} and consider $K=\langle
  \kappa\rangle$ and $g=\rho$. Again, this produces a much higer
  dimensional space than the original example in $\R^3$.
\end{example}

\begin{example}
  The homoclinic triples for wreath products defined in Example
  \ref{eg:wreath} with $\mathcal{L}$ finite satisfy the hypothesis
  [HTL] of Theorem \ref{thmsufficient}.  There are therefore open sets
  of equivariant vector fields on the complex regular representation
  with \rht s for which these triples occur.
\end{example}

\begin{remark} 
  Many examples of robust \RHT s in the literature are for irreducible
  representations (for example, the Guckenheimer-Holmes example in
  Section~\ref{secGHmodel}) or for group actions where there are only
  a few irreducible components (for example, the Kuramoto-Sivashinsky
  example in Section~\ref{secKSmodel} where there are two).  Our
  result, Theorem~\ref{thmsufffinite}, by contrast gives robust \RHT s
  on much larger dimensional spaces with many irreducible components.
  Usually, no irreducible representation will be sufficiently rich to
  satisfy the hypothesis of Theorem \ref{thmsufficient}, and this begs
  the question of how to obtain optimally small estimates for the
  dimension of the group action one needs to consider obtain robust
  \RHT s for a given group.
\end{remark}

\begin{remark}
  If $G$ is a finite group that satisfies Theorem~\ref{thmsufffinite}
  we say that $G$ is robust. Clearly any finite group that has a
  robust subgroup is robust, giving another way of constructing many
  groups admitting \RHT s.
\end{remark}

%%%%%%%%%%
\subsection{An example with $G=\D_8$} 
\label{sec:D_8-example}

We present an example demonstrating the difference between the
necessary conditions of Proposition \ref{proptwist} and the sufficient
conditions of Theorem \ref{thmsufficient}.  It is in some sense
the simplest example
of a homoclinic triple that does not satisfy the extra sufficient
condition [HTL] of Theorem \ref{thmsufficient}, and appears in the
$\D_8$-table in \S\ref{sec:dihedral}. We consider complex
representations so that the global connectivity hypothesis [HTG] of
Theorem \ref{thmsufficient} is automatically satisfied, as in Theorem
\ref{thmsufffinite}.

Using the notation for the dihedral groups introduced in
\S\ref{sec:dihedral}, consider the homoclinic triple $(K,g,H)$, where
$K=\langle\kappa\rangle$, $g=\rho$ and $H=\langle\kappa,
\rho^2\rangle$.  This is a homoclinic triple by Theorem \ref{thm:D4n},
and indeed it is a minimal triple.

Consider the irreducible representations (irreps) of $G=\D_8$, of
which there are 7 in all (Table \ref{table:D_8}).  It can be seen from
the table that, for a representation $V$ of $\D_8$, $H$ is an isotropy
subgroup if and only if $V$ contains at least one copy of $A_2$, and
$K$ is an isotropy subgroup if and only if $V$ contains at least one
copy of either $E_1$ or $E_2$.

\begin{table}[ht]
  $$\begin{array}{c|ccccccc||c||c}
    & e & (\kappa) & (\kappa\rho) & (\rho) & (\rho^2) & (\rho^3) 
    & \rho^4 & \D_4 & \Z_2\cr
    \hline
    A_0 & 1 & 1 & 1 & 1 &1 & 1 & 1 & A_0' & A_0'' \cr
    A_1 & 1 & -1 & -1 & 1 &1 & 1 & 1 & A_1' & A_1'' \cr 
    A_2 & 1 & 1 & -1 & -1 &1 & -1 & 1 & A_0' & A_0'' \cr 
    A_3 & 1 & -1 & 1 & -1 &1 & -1 & 1 & A_1' & A_1'' \cr 
    \hline
    E_1 & 2 & 0 & 0 & \sqrt{2} & 0 & -\sqrt{2} & 0 & E' & A_0'' + A_1''\cr
    E_2 & 2 & 0 & 0 & -\sqrt{2} & 0 & \sqrt{2} & 0 & E' & A_0'' + A_1'' \cr
    F & 2 & 0 & 0 & 0 & -2 & 0 & 2 & A_2'+A_3' & A_0'' + A_1'' \cr
    \hline
  \end{array}$$
  
  \caption{Character table for $\D_8$.  The final two columns give the
    corresponding representations of $\D_4=\langle\kappa,\rho^2\rangle$ 
    and $\Z_2 = \langle\kappa\rangle$ via restriction (the characters 
    for the $\D_4$ and $\Z_2$ representations can be deduced from this table).}
  
  \label{table:D_8} 
\end{table}

\begin{theorem} 
  Let $V$ be a complex representation of $\D_8$ for which $K$ and $H$
  are isotropy subgroups, where $K=\Z_2$ and $H=\D_4$ as above.  There
  is an open set of equivariant vector fields with homoclinic triple
  $(K,g,H)$ with $g=\rho$ if and only if $V$ contains at least one
  copy of the irreducible representation $F$ (see Table
  \ref{table:D_8}).
\end{theorem}

\begin{proof} 
  This is based on the proof of Theorem \ref{thmsufficient} and Remark
  \ref{rmk:HTL}.  Let $x\in \Fix(H)$. Then $T_xV=V$ and the
  linearization $L$ of the equivariant vector field at $x$ is
  $H$-equivariant, and we need to consider the eigenvalues of $L$ on
  each irrep. Consider the isotypic decomposition of $V$:
  $$V = a_0A_0 \oplus a_1 A_1 \oplus a_2 A_2 \oplus a_3 A_3 \oplus 
  e_1E_1 \oplus e_2E_2 \oplus fF,$$%
  where $a_0,\dots,f$ are non-negative integers.  Then
  $$\begin{array}{rcl}
    \Fix(H,V) &=& a_0A_0 \oplus a_2A_2 \\
    \Fix(K,V) &=& a_0A_0 \oplus a_2A_2 \oplus e_1E_1^{(s)} 
         \oplus e_2E_2^{(s)} \oplus fA_2' \\
    \Fix(\Kginv,V) &=& a_0A_0 \oplus a_2A_2 \oplus e_1E_1^{(d)} 
         \oplus e_2E_2^{(d)} \oplus fA_3'.
  \end{array}$$
  Some of this notation needs explaining.  The $A_k$ parts should be
  self-explanatory. The two representations $E_1$ and $E_2$ of $\D_8$
  are the usual symmetry groups of the regular octagon; in the first
  $\rho$ acts by rotation by $\pi/4$ while in $E_2$ it acts by
  rotation by $3\pi/4$. Restricting the action to $H$, picks out a
  square in the octagon, whose vertices lie at alternating vertices of
  the octagon say, and then $E_j^{(s)}$ is a line of reflection in
  $E_j$ passing through midpoints of 
  a pair of sides of the square, while
  $E_j^{(d)}$ is a diagonal line of the square.  Note that any
  $H$-equivariant linear vector field on $E_j'$ has the same
  eigenvalues on both $E_j^{(s)}$ and $E_j^{(d)}$, so if $f=0$ in the
  representation, the eigenvalues on $\Fix(\Kginv,V)$ and $\Fix(K,V)$
  cannot be distinct.
  
  On the other hand, the irrep $F$ decomposes into two $H$-irreps,
  $F=A_2'+A_3'$. Indeed, in the action of $\D_8$ on $F$, $\rho$ acts
  as rotation by $\pi/2$, so that $\D_8$ acts as the symmetry group of
  the square, and the action of $H$ is just by a pair of reflections.
  These two reflections are in fact those in $K$ and $\Kginv$, and so
  $\Fix(K,F) = A_2'$ say, and $\Fix(\Kginv,F) = A_3'$.  It follows
  that an $H$-equivariant $L$ can be chosen so that the eigenvalues on
  $A_2'$ and $A_3'$ are of opposite sign, and there is an open set of
  such $L$.  Consequently if $f\geq 1$ then there is an open
  set of linear vector fields $L$ satisfying (\ref{eq:dimensions}), as
  required by Remark \ref{rmk:HTL}.
\end{proof}

%%%%%%%%%%%%%%%%%%%%%%%%%%
\section{Discussion}
\label{secdiscuss} \setcounter{equation}{0}

There has been much work on robust heteroclinic cycles that has looked
at structure and existence of \RHT s for given group representations.
What we have attempted here is to understand better the
group-theoretic conditions on a group necessary to find a
representation admitting \RHT s. We have found necessary conditions in
Proposition~\ref{proptwist} and sufficient conditions in
Theorem~\ref{thmsufficient} but there remains a gap in the hypotheses
that it would be nice to be able to close. Specifically one would like
to be able to characterise a weaker version of [HTL] that would be
both necessary and sufficient, and Remark \ref{rmk:HTL} together with
the example in \S\ref{sec:D_8-example} shows that this condition must
include some information on the local structure of the action.

A number of other questions are suggested by this study. As mentioned
already, there are optimality questions; for example,
given a group with a homoclinic triple that does admit robust \RHT s, 
how small a representation can one 
consider to find a robust \RHT? Also, how rare are robust \RHT s for
equivariant systems. For example, what proportion of finite groups of
order $n$ have the necessary complexity to admit robust \RHT s? What
proportion of irreducible representations of finite groups of
dimension $n$ admit robust \RHT s?

It does not seem to be a trivial task to extend or generalize the
results here to apply to robust heteroclinic cycles due partly to the
fact that the interconnection possibilities are much greater and
ensuring their robustness requires that many local and global
conditions are fulfilled simultaneously (see for example \cite{MCG89}).
However, it should be possible to generalize the results to apply to
heteroclinic cycles between classes of more general transitive
invariant sets.

%%%%%%%%%%%%%%%%%%%%%%%%%%%%%%%%

\end{document}